\documentclass[a4paper]{article}

\usepackage{amsmath,amsthm,amsfonts}

\usepackage{cite}

\theoremstyle{plain}

\newtheorem{prop}{Proposition}
\newtheorem{lem}[prop]{Lemma}
\newtheorem{thm}[prop]{Theorem}

\theoremstyle{remark}

\newtheorem{rem}[prop]{\bf Remark}

\newcommand{\CC}{\mathbb{C}}

\newcommand{\RR}{\mathbb{R}}

\newcommand{\HH}{\mathcal{H}}

\newcommand{\GG}{\mathcal{G}}

\newcommand{\LL}{\mathcal{L}}

\newcommand{\NN}{\mathcal{N}}

\newcommand{\ZZ}{\mathbb{Z}}

\newcommand{\ran}{\mathop{\mathrm{ran}}}

\newcommand{\dom}{\mathop{\mathrm{dom}}}

\newcommand{\res}{\mathop{\mathrm{res}}}

\newcommand{\spec}{\mathop{\mathrm{spec}}}

\newcommand{\gr}{\mathop{\mathrm{gr\,}}}

\newcommand{\diag}{\mathop{\mathrm{diag\,}}}

\newcommand{\sgn}{\mathop{\mathrm{sgn\,}}}

\title{\bf Resolvents of self-adjoint extensions
with mixed boundary conditions}

\author{Konstantin Pankrashkin\\[2ex]Institut f\"ur Mathematik, Humboldt-Universit\"at zu Berlin\\
Rudower Chaussee~25, 12489 Berlin, Germany\\
E-mail: const@mathematik.hu-berlin.de}

\begin{document}


\date{}

\maketitle

\begin{abstract} We prove a variant of Krein's resolvent formula
expressing the resolvents of self-adjoint extensions
through the associated boundary conditions.
Applications to solvable quantum-mechanical problems are discussed.
\end{abstract}



\section{Introduction}

Self-adjoint extensions of symmetric operators arise in many
areas of mathematical physics, like solvable models of quantum mechanics
and quantization problems~\cite{pavlov,alb,ak,gtv}.
Of interest are first of all their spectral and scattering properties. 
One of the most powerful tools
for such an analysis is the well-known Krein's formula connecting
the resolvents of two self-adjoint extensions of a symmetric operator, see e.g. \cite{ag,dm,GMT,krein,LT,BKN,pos1}.
It is well known that all self-adjoint extensions of a densely
defined symmetric operator $S$
with equal deficiency indices are parameterized by unitary operators in a certain auxiliary
Hilbert space. In many situations there is a distinguished self-adjoint extension $H_0$
whose spectral properties are known (for example, the free Laplacian
the theory of zero-range potentials), and it would be useful to analyze
all other self-adjoint extensions in terms of $H_0$ and
of the unitary parameters. From this point of view, the existing versions of the resolvent
formula are either applicable to disjoint extensions only (i.e. with special unitary parameters) or involve multivalued
operators, and in both cases the global contribution of the parameters
to the spectral properties of the extensions remains unclear.
On the other hand, only considering a complete family of self-adjoint extensions
delivers a relevant description of topology change phenomena
in quantum-mechanical problems~\cite{AIM,CFT}. We provide a resolvent
formula covering the whole family of self-adjoint extensions
in an explicit form in the present report.

We look at the situation from an abstract point of view.
Let $S$ be a closed densely defined symmetric operator with the
deficiency indices $(n,n)$, $0<n\le\infty$, acting on a 
Hilbert space $\HH$. Let $\GG$ be an auxiliary Hilbert space
such that $\dim \GG=n$. One says that a triple
$(\GG,\Gamma_1,\Gamma_2)$, where $\Gamma_1$ and
$\Gamma_2$ are linear maps from the domain $\dom S^*$ of the adjoint of $S$
to $\GG$, is a \emph{boundary triple} (or a \emph{boundary value space})
for $S$ if the following two conditions are satisfied:
\begin{itemize}
\item for any $\phi,\psi\in\dom S^*$ there holds an abstract integration
by parts,
\[\langle \phi,
S^*\psi\rangle- \langle S^*\phi, \psi\rangle= \langle \Gamma_1
\phi,\Gamma_2\psi\rangle-\langle \Gamma_2
\phi,\Gamma_1\psi\rangle,\] 
\item the map $(\Gamma_1,\Gamma_2):\dom S^*\to \GG\oplus \GG$ is surjective.
\end{itemize}
(We assume that the inner products in $\HH$ and $\GG$
are linear with respect to the second argument.)
It has been known for a long time that all self-adjoint extensions of $S$ are parameterized
by self-adjoint linear relations in $\GG\oplus\GG$ \cite{dm,gorb}.
On the other hand, in many situations
it is natural to parameterize self-adjoint extensions by boundary conditions
of the form $A\Gamma_1\phi=B\Gamma_2\phi$, where $A$ and $B$
are bounded linear operators acting on $\GG$ (an important class of such problems
comes from the study of quantum graphs an hybrid manifolds~\cite{ks,bg,E05,bgl,Ku}). 
In the present note we show that the resolvent formula
admits a simple form in terms of $A$ and $B$, which permits
to cover the whole family of self-adjoint extensions of $S$ (see theorems \ref{prop-main}
and \ref{prop-main2}).
The case when $S$ has finite deficiency indices was considered in~\cite{ap}, and here we are interested in the infinite dimensional situation.
In section~\ref{sec4} we discuss some quantum-mechanical examples in which
mixed boundary conditions arise.

\section{Parameterization of self-adjoint linear relations}

Let us recall some basic facts on linear relations. For a more detailed
discussion we refer to~\cite{arens}.
Any linear subspace of $\GG\oplus\GG$ will be called a \emph{linear relation} on $\GG$.
For a linear relation $\Lambda$ on $\GG$ the sets
\begin{align*}
\dom\Lambda&=\{x\in \GG:\,\exists y\in \GG
\text{ with } (x,y)\in\Lambda)\},\\
\ran \Lambda&=\{
x\in \GG:\,\exists y\in \GG \text{ with } (y,x)\in\Lambda) \},\\
\ker\Lambda&=\{x\in \GG: (x,0)\in\Lambda\}
\end{align*}
will be called the \emph{domain}, the \emph{range}, and the \emph{kernel} of~$\Lambda$,
respectively.  The linear relations
\begin{align*}
\Lambda^{-1}&=\{ (x,y):\,(y,x)\in\Lambda\},\\
\Lambda^*&=\{(x_1,x_2): \langle x_1,y_2\rangle=\langle x_2,y_1\rangle
\quad\forall (y_1,y_2)\in\Lambda\}
\end{align*}
are called \emph{inverse} and \emph{adjoint} to $\Lambda$, respectively.
For $\alpha\in\CC$ we put $\alpha \Lambda=\{(x,\alpha
y):\,(x,y)\in\Lambda\}$.
For two linear relations
$\Lambda',\Lambda''\subset \GG\oplus \GG$ one can define their \emph{sum}
\[
\Lambda'+\Lambda''=\{(x,y'+y''),\,(x,y')\in\Lambda',\,
(x,y'')\in\Lambda''\};
\]
clearly, one has $\dom
(\Lambda'+\Lambda'')=\dom\Lambda'\cap\dom\Lambda''$. The graph of
any linear operator $L$ on $\GG$ is a linear relation, which
we denote by $\gr L$. Clearly, if $L$ is invertible,
then $\gr L^{-1}=(\gr L)^{-1}$. For arbitrary linear operators
$L',L''$ one has $\gr L'+\gr L''=\gr (L'+L'')$. Therefore, the set
of linear operators is naturally embedded into the set of linear
relations. 
In what follows we consider mostly only closed linear relations, i.e. which are closed
linear subspaces in $\GG\oplus\GG$. Clearly, this notion generalizes the notion
of a closed operator. By analogy with operators, one introduces
the notion of the \emph{resolvent set} $\res\Lambda$ of a
closed linear relation $\Lambda$ by the rule
$\res\Lambda=\{\lambda\in\CC: \ker(\Lambda-\lambda I)=0
\text{\,\, and } \ran(\Lambda-\lambda I)=\GG\}$,
where $I\equiv \gr \mathop{\mathrm{id}_\GG}=\big\{(x,x),\, x\in\GG\big\}$.
In other words, the condition $\lambda\in\res\Lambda$
means that $(\Lambda-\lambda I)^{-1}$ is the graph of
a certain linear operator
defined everywhere; this operator is bounded
due to the closed graph theorem.

A linear relation $\Lambda$ on $\GG$ is called
 \emph{symmetric} if $\Lambda\subset
\Lambda^*$ and is called \emph{self-adjoint} if
$\Lambda=\Lambda^*$ (in the geometric language, they are
called \emph{isotropic} and \emph{Lagrangian subspaces}, respectively,
see remark~\ref{rem8} below).
A linear operator $L$ in $\GG$ is
symmetric (respectively, self-adjoint), iff its graph is a
symmetric (respectively, self-adjoint) linear relation.
A self-adjoint linear relation (abbreviated as s.a.l.r.)
is always maximal symmetric, but the converse in not true;
examples are given by the graphs of maximal symmetric operators with 
deficiency indices $(m,0)$, $m>0$.

Our aim now is to find a suitable way for presenting
s.a.l.r. Let $A$, $B$ be bounded linear operators on $\GG$.
We introduce the notation
\[
\Lambda^{A,B}=\big\{ (x_1,x_2)\in \GG\oplus \GG, \quad  A x_1=B x_2
\big\}.
\]
We say that a linear relation $\Lambda$ \emph{is parameterized by
the operators $A$ and $B$} if $\Lambda=\Lambda^{A,B}$.
Conditions for $\Lambda^{A,B}$ to be self-adjoint can be written
is many ways, see e.g. \cite{dm2,RB}. We will
use the conditions obtained in~\cite{bg}.
\begin{prop}[Proposition B in~\cite{bg}]
Denote by $M^{A,B}$
an operator acting on $\GG\oplus\GG$ by the rule
\begin{equation}
        \label{eq-mab}
M^{A,B}=\begin{pmatrix}
A & -B\\
B &A
\end{pmatrix},
\end{equation}
then
the linear relation $\Lambda^{A,B}$ is
self-adjoint iff $A$ and $B$ satisfy the following two conditions:
\begin{gather}
        \label{eq-bg1}
	AB^*=BA^*,\\
        \label{eq-bg2}
	\ker M^{A,B}=0.
\end{gather}
\end{prop}
\begin{prop}[Theorem~3.1.4 in~\cite{gorb}]
           \label{prop-iu}
For a given linear relation $\Lambda$ in $\GG$ there is
a unique unitary operator $U$ in $\GG$ (called the \emph{Cayley transform}
of $\Lambda$) such that
the condition $(x_1,x_2)\in\Lambda$ is equivalent to
$i(1+U)x_1=(1-U)x_2$, i.e. $\Lambda=\Lambda^{i(1+U),1-U}$.
\end{prop}
Taking $U$ in the form $U=e^{-2i\Phi}$, where $\Phi$ is a self-adjoint operator
in $\GG$, one write any s.a.l.r. as $\Lambda^{\cos\Phi,\sin\Phi}$.

Although proposition~\ref{prop-iu} claims that
there exists a one-to-one correspondence between
s.a.l.r.s and unitary operators, for a given
s.a.l.r. $\Lambda$ it is difficult to find its
Cayley transform,
but there are many other ways to represent it as $\Lambda^{A,B}$
with suitable $A$ and $B$.

In what follows we will need a parameterization of s.a.l.r. satisfying stronger
conditions than \eqref{eq-bg1} and~\eqref{eq-bg2}. More precisely, we replace
the condition~\eqref{eq-bg2} by
\begin{equation} \label{eq-bg3} 0\in \res
M^{A,B}. \end{equation}
We say that a pair of bounded operators $A$ and $B$ satisfying~\eqref{eq-bg1}
and~\eqref{eq-bg2} is \emph{normalized} if the condition~\eqref{eq-bg3} is satisfied.
Clearly, in the case of finite-dimensional $\GG$ the
conditions~\eqref{eq-bg2} and~\eqref{eq-bg3} are equivalent. Moreover, in this
case these conditions are equivalent to the following one~\cite{ks}:
\[ \text{the $n\times 2n$ matrix $(AB)$ has maximal rank.}
\]
(Note that this can be written also as $\det(A A^*+B B^*)\ne0$,
which can be found in the textbooks on operator theory~\cite[Section~125, Theorem~4]{ag}).
In general, the conditions \eqref{eq-bg2} and~\eqref{eq-bg3} do not coincide:
if one replaces $A$ by $LA$ and $B$ by $LB$,
where $L$ is a bounded linear operators with
$\ker L=0$ and $(\ran L)^\bot\ne 0$, this does not change
the subspace $\Lambda^{A,B}$, but the condition~\eqref{eq-bg3} will
not be satisfied. Moreover, this construction is the only source of
``denormalization''.

\begin{prop}\label{prop-lab}
Let $A$, $B$, $C$, $D$ be bounded operators in $\GG$ and
$\Lambda$ be a s.a.l.r in $\GG$ such that $\Lambda=\Lambda^{A,B}=\Lambda^{C,D}$.
Assume that $A$ and $B$ are normalized, then
there exists a bounded injective operator $L$ on $\HH$
with $C=LA$ and $C=LB$.
\end{prop}

\begin{proof}
Introduce operators $M_1,M_2:\GG\oplus\GG\to\GG$ by
$M_1(x_1,x_2)=Ax_1-Bx_2$, $M_2(x_1,x_2)=Cx_1-Dx_2$, $x_1,x_2\in\GG$.
The condition~\eqref{eq-bg3} says, in particular, that
$\ran M_1=\GG$.
Clearly, the null spaces of $M_1$ and $M_2$ coincide,
$\ker M_1=\ker M_2=\Lambda$, therefore, there exists an injective
operator $L$ in $\HH$ such that $M_2=LM_1$, which implies
$C=LA$ and $D=LB$. Let us show that $L$
is bounded. Use the notation~\eqref{eq-mab}, then
there holds $M^{C,D}=(L\oplus L)M^{A,B}$.
Due to the condition~\eqref{eq-bg3} the operator $M^{A,B}$
has a bounded inverse defined everywhere. Therefore,
$L\oplus L=M^{C,D}\big(M^{A,B}\big)^{-1}$ is bounded
and defined everywhere, so is $L$.
\end{proof}
It is important to emphasize that for a given s.a.l.r one can always find
a normalized parameterization, as the following proposition shows.
\begin{prop}\label{prop-AB}
\textup{(a)} Let $U$ be a unitary operator in $\GG$,
then the operators $A=i(1+U)$ and $B=1-U$ satisfy the conditions
\eqref{eq-bg1} and~\eqref{eq-bg3}.

\textup{(b)} Any s.a.l.r. can be parameterized by operators
$A$ and $B$ satisfying~\eqref{eq-bg1} and~\eqref{eq-bg3}.
\end{prop}

\begin{proof}
(a)
The condition~\eqref{eq-bg1} is obviously satisfied, so we prove only~\eqref{eq-bg3}.
First of all note that the operator $M^*$ adjoint to
$M=M^{A,B}$ is given by the following operator-matrix:
\[
M^*=\begin{pmatrix}
A^* & B^*\\
-B^* & A^*,
\end{pmatrix}
\]
or, in our case,
\[
M^*=\begin{pmatrix}
-i(1+U^*) & 1-U^*\\
U^*-1 & -i(1+U^*)
\end{pmatrix}.
\]
Let us show that $\ker M^*=0$. Assume $x=(x_1,x_2)\in\ker M^*$, $x_1,x_2\in\GG$,
then
\begin{align}
         \label{u1}
-i(1+U^*)x_1+(1-U^*)x_2&=0,\\
         \label{u2}
(U^*-1)x_1 -i(1+U^*)x_2&=0.
\end{align}
Multiplying~\eqref{u1} by $i$ and adding the result to~\eqref{u2}
one arrives at $U^*(x_1-ix_2)=0$; as $U^*$ is unitary, we have $x_1-ix_2=0$.
On the other hand, multiplying~\eqref{u1} by $i$ again and subtracting~\eqref{u2}
from it, we obtain $x+i x_2=0$, which says that $x_1=x_2=0$.

As $(\ran M)^\bot=\ker M^*$, the linear subspace
$\ran M$ is dense in $\GG$. Now to prove \eqref{eq-bg3}
it is sufficient to show that for any sequence
$(x^n)\in \GG\oplus\GG$, $x^n=(x_1^n,x_2^n)$,
$x_1^n,x_2^n\in\GG$, the condition $\lim_{n\to\infty}M x^n= 0$
implies the convergence of $(x^n)$ to $0$, which we will do now.

Assuming the existence of the limits
\[
\lim_{n\to\infty}\big(i(1+U)x_1^n+(U-1)x_2^n\big)= 0,\quad
\lim_{n\to\infty}\big((1-U)x_1^n+i(1+U)x_2^n\big)=0
\]
one sees immediately that the sequences
$(x^n_1+i x^n_2)$ and $\big(U(-x^n_1+ix^n_2)\big)$ converge to $0$.
As $U$ is unitary,
the sequence $(-x^n_1+i x^n_2)$ converges to $0$ too,
which shows that $\lim_{n\to\infty} x^n_1=\lim_{n\to\infty}x^n_2=0$.

(b) This is an obvious corollary of (a) and proposition~\ref{prop-iu}.
\end{proof}

Finally, we are able to give another description of a s.a.l.r.
with the help of its normalized parameterization.

\begin{lem} \label{lem-LAB}
Let bounded operators $A$, $B$ parameterize a s.a.l.r. in $\GG$
and be normalized, then $\Lambda^{A,B}=\big\{(B^*u,A^*u),\,u\in\GG\big\}$.
\end{lem}
\begin{proof}
Set $\Lambda':=\big\{(B^*u,A^*u),\,u\in\GG\}$. Clearly, due to \eqref{eq-bg1}
there holds the inclusion $\Lambda'\subset\Lambda^{A,B}$.
Let us show that $\Lambda^{A,B}=\Lambda'$.
The condition~\eqref{eq-bg3} means,
in particular, that the operator $(M^{A,B})^*$ has a bounded inverse
and, therefore, maps closed sets to closed sets.
As $\Lambda'=(M^{A,B})^*(0\oplus \GG)$, $\Lambda'$ is closed.
As $\Lambda^{A,B}$ is also closed, it is sufficient to prove
that $(\Lambda')^\bot\cap\Lambda^{A,B}=0$. Assume $x=(x_1,x_2)\in
(\Lambda')^\bot\cap\Lambda^{A,B}$, $x_1,x_2\in\GG$. The condition
$x\in\Lambda^{A,B}$ means that $Ax_1-Bx_2=0$, and the equality
$\langle x,y\rangle=0$ for any $y\in \Lambda'$ results
in $\langle x_1,B^*u\rangle+\langle x_2,A^*u\rangle=0$ 
or $\langle Bx_1+Ax_2,u\rangle=0$ for any $u\in\GG$, i.e.
$B x_1+Ax_2=0$. Therefore, $M^{A,B}x=0$ and due to~\eqref{eq-bg3}
there holds $x=0$.
\end{proof}

\section{Resolvents of self-adjoint extensions}

The language of linear relations is widely used in the theory of
self-adjoint extensions of symmetric operators~\cite{RB,gorb,Kc}.
We point out that any symmetric operator with equal deficiency indices
(finite or infinite) has a boundary triple~\cite[Theorem 3.1.5]{gorb}.

\begin{prop}[Theorem 3.1.6 in \cite{gorb}]
Let $S$ be a closed symmetric operator with equal deficiency indices
acting on a certain Hilbert space,
and $(\GG,\Gamma_1,\Gamma_2)$ be its boundary triple, then
there
is a bijection between all self-adjoint extensions of $S$ and
s.a.l.r's on $\GG$. A self-adjoint extension $H^\Lambda$
corresponding to a s.a.l.r. $\Lambda$ is the restriction of $S^*$ to
elements $\phi\in\dom S^*$ satisfying the abstract boundary conditions
$(\Gamma_1\phi,\Gamma_2\phi)\in\Lambda$.
\end{prop}

To investigate spectral properties of the self-adjoint extensions
it is useful to know their resolvents.
To write Krein's formula for the resolvents we need some
additional constructions~\cite{dm}. For $z\in\CC\setminus\RR$, let $\NN_z$
denote the corresponding deficiency subspace for $S$, i.e.
$\NN_z=\ker(S^*-z)$. The restrictions of $\Gamma_1$ and $\Gamma_2$
onto $\NN_z$ are invertible linear maps from $\NN_z$ to $\GG$. Put
$\gamma(z)=\big(\Gamma_1|_{\NN_z}\big)^{-1}$ and
$Q(z)=\Gamma_2\gamma(z)$; these maps form holomorphic families from
$\CC\setminus\RR$ to the spaces $\LL(\GG,\HH)$ and $\LL(\GG,\GG)$ of
bounded linear operators from $\GG$ to $\HH$ and from $\GG$ to $\GG$ respectively.
Denote by $H^0$ the self-adjoint extension of $S$ given by the
boundary condition $\Gamma_1\phi=0$,
then the maps $\gamma(z)$ and $Q(z)$ have analytic continuations to
the resolvent set $\res H^0$, and for all $z,\zeta\in\res H^0$ one
has, in particular,
\begin{equation}
        \label{q-fun}
Q(z)-Q^*(\zeta)=(z-\overline\zeta)\,\gamma^{\,*}(\zeta)\,\gamma(z).
\end{equation}
The maps $\gamma(z)$  and $Q(z)$ are called the \emph{$\Gamma$-field} and
the \emph{$\mathcal{Q}$-function} for the pair $(S,H^0)$, respectively~\cite{GMT,krein,LT}.
(The $\mathcal{Q}$-function is called sometimes the \emph{Weyl $M$-function}
of the boundary triple $(\GG,\Gamma_1,\Gamma_2)$~\cite{dm,ABMN,BMN}.)
Similar objects arise naturally also in the study of singular perturbations
of self-adjoint operators~\cite{pos1,pos2}.

The following proposition describes the resolvents of
the self-adjoint extensions of $S$.
\begin{prop}[Krein's resolvent formula, cf. Propositions~1 and~2 in~\cite{dm}]\label{prop-dm}
Let $H^\Lambda$ be a self-adjoint extension of $S$,
which is the restriction of $S^*$ to the set of functions $\phi\in\dom S^*$
satisfying $(\Gamma_1\phi,\Gamma_2\phi)\in\Lambda$,
where $\Lambda$ is a s.a.l.r. in $\GG$. Then a number
$z\in \res H^0$ lies in the spectrum
of $H^\Lambda$ iff $0\notin \res\big(\gr Q(z)-\Lambda\big)^{-1}$.
For any $z\in\res H^0\cap\res H^\Lambda$ there holds
\begin{equation}
     \label{krein}
(H^\Lambda-z)^{-1}=(H^0-z)^{-1}-\gamma(z)\, C_\Lambda(z)\,
\gamma^{\,*}(\bar{z}),
\end{equation}
where $C_\Lambda(z)$ is a bounded linear operator on $\GG$
with $\gr C_{\Lambda}(z)=\big(\gr Q(z)-\Lambda\big)^{-1}$.
\end{prop}
It is worth emphasizing that the correspondence between the spectral types
of $H^\Lambda$ and $C_\Lambda(z)$ (discrete spectra, essential spectra etc.)
is a rather difficult problem, cf. \cite{BMN,ABMN,Br}

The calculation of $C_\Lambda(z)$ is a rather difficult technical
problem, as it involves ``generalized'' operations with linear
relations. Such difficulties do not arise if $\Lambda$ is the graph of
a certain self-adjoint linear operator $L$ (i.e. if $\Lambda$ can be
injectively projected onto $\GG\oplus 0$);
the boundary conditions take the form
$\Gamma_2\phi=L\,\Gamma_1\phi$,
and such extensions are called \emph{disjoint to
$H^0$} because of the equality $\dom H^\Lambda\cap\dom
H^0=\dom S$ (the operator $S$ is then called the
\emph{maximal common part} of $H^0$ and $H^\Lambda$). Then the subspace
$\gr Q(z)-\Lambda$ is the graph of the  invertible \emph{operator}
$Q(z)-L$, and $C^\Lambda(z)=\big(Q(z)-L\big)^{-1}$.

\begin{rem}\label{rem8}
One uses sometimes a different terminology, which is more related to geometry.
The space $\GG\oplus\GG$ is equipped with a symplectic structure
given by the skew-linear form $[\cdot,\cdot]$,
$[(p_1, q_1),(p_2,q_2)]=\langle p_1,q_2\rangle
-\langle p_2, q_1\rangle$, $p_1,p_2,q_1,q_2\in\GG$.
In these terms, symmetric linear relations are linear subspaces on which
this form vanishes (they are called more often
\emph{isotropic subspaces}) and s.a.l.r.s are \emph{Lagrangian subspaces}
(i.e. those coinciding with their skew-orthogonal complements
with respect to the form $[\cdot,\cdot]$).
In the case of real finite-dimensional $\GG$, such objects
appeared in the semiclassical analysis~\cite{mf}.
They play an important role in the
description of classical dynamics, as invariant manifolds of
integrable Hamiltonian systems are Lagrangian, i.e. all tangent spaces are
Lagrangian. If $e_j$, $j=1,\dots,n$, form an orthogonal basis
in $\GG$, then the $2n$ vectors $(e_j,0)$, $(0,e_j)$, $j=1,\dots,n$,
form a symplectic basis in $\GG\oplus\GG$. Let $\theta$ be a subset of $\{1,\dots,n\}$,
then the linear hull of the vectors $(e_j,0)$, $j\notin\theta$, $(0,e_j)$, $j\in\theta$,
is called a \emph{coordinate subspace}.
Arnold's lemma~\cite{arn} says that an arbitrary
Lagrangian subspace can be injectively projected onto one of the coordinate subspaces.
This is one of the central points in the construction of WKB-solutions resulting
in the Bohr-Sommerfeld-Maslov quantization rule~\cite{mf}. In order to
reduce calculations, one usually tries to minimize the number of elements in $\theta$,
as this number is, roughly speaking, the number of partial Fourier transforms
needed to write a formula for the solution.

Arnold's lemma can be transferred to the case of
complex finite-dimensional $\GG$ and applies
to symmetric operators with equal and finite deficiency indices
as follows. Start with an arbitrary boundary value space $(\GG,\Gamma_1,\Gamma_2)$.
Fix an orthogonal basis $(e_j)$ in $\GG$ and denote by $\Gamma^{j}_k$
the $j$th component of $\Gamma_k$ in this basis,
$\Gamma_k^j=\langle e_j,\Gamma_k\,\cdot\rangle$,
$j=1,\dots n$, $k=1,2$. For a subset $\theta\in\{1,\dots,n\}$
we define new boundary operators $(\Gamma_1^\theta,\Gamma_2^\theta)$ by
\[
(\Gamma_1^\theta)^j=\Gamma_1^j,\quad (\Gamma_2^\theta)^j=\Gamma_2^j, \quad j\notin\theta,
\quad (\Gamma_1^\theta)^j=\Gamma_2^j,\quad (\Gamma_2^\theta)^j=-\Gamma_1^j,\quad j\in\theta.
\]
Clearly, the triple $(\GG,\Gamma_1^\theta,\Gamma_2^\theta)$ is a new boundary value space,
and for a \emph{fixed} self-adjoint extension $H^\Lambda$ one can choose $\theta$
for which the boundary conditions for $H^\Lambda$ take the form
$\Gamma^\theta_2\phi=L\Gamma_1^\theta\phi$ with a certain matrix $L$.
For each of these boundary value spaces $(\GG,\Gamma_1^\theta,\Gamma_2^\theta)$
one should recalculate the maps $\gamma(z)$ and $Q(z)$ entering the resolvent
formula, which brings a number of calculations, cf. \cite{ABMN,GMT}.
\end{rem}

As we have shown in proposition~\ref{prop-AB}, all self-adjoint boundary
conditions can be represented with the help of two bounded linear operators
$A$ and $B$ by
\begin{equation}
     \label{A-B}
A\Gamma_1\phi=B\Gamma_2\phi \quad \Leftrightarrow \quad
(\Gamma_1\phi,\Gamma_2\phi)\in\Lambda^{A,B},
\end{equation}
where $A$ and $B$ satisfy~\eqref{eq-bg1} and~\eqref{eq-bg3}. Our
aim is to show that the resolvent formula~\eqref{krein} admits a simple form in
terms of these two operators.

\begin{lem} \label{lem-main}
 Let bounded operators $A$, $B$ parameterize
a s.a.l.r. in $\GG$ and be normalized. Then for any $z\in\res H^0$
the following three conditions are equivalent:
\begin{itemize}
\item[\rm (a)] $0\in\res\big(\gr Q(z)-\Lambda^{A,B}\big)$,
\item[\rm (b)] $0\in\res \big(B Q(z)-A\big)$,
\item[\rm (c)] $0\in\res \big(Q(z)B^*-A^*\big)$.
\end{itemize}
If these conditions are satisfied, then
\begin{equation}
          \label{eq-QLAB}
\big(\gr Q(z)-\Lambda^{A,B}\big)^{-1}= \gr B^*\big(Q(z)B^*-A^*\big)^{-1}=\gr\big(BQ(z)-A\big)^{-1}B.
\end{equation}
\end{lem}

\begin{proof}
Let us express the linear relation $\gr Q(z)-\Lambda^{A,B}$
through $A$ and $B$. Due to lemma~\ref{lem-LAB} one has
$\Lambda^{A,B}=\{(B^*u,A^*u),\,u\in\GG\}$. Therefore,
$\dom \big(\gr Q(z)-\Lambda^{A,B}\big)=\ran B^*$, and
there holds
\begin{equation}
       \label{eq-qab}
\big(\gr Q(z)-\Lambda^{A,B}\big)
=\Big\{\big(B^*u,\, Q(z)B^*u-A^*u\big),\,\,u\in\GG\Big\}.
\end{equation}

Assume that (a) is satisfied and show (c).
Clearly, $\ran (\gr Q(z)-\Lambda^{A,B})=\ran \big(Q(z)B^*-A^*\big)$,
and there holds $\ran \big(Q(z)B^*-A^*\big)=\GG$.
We show now that $\ker \big(Q(z)B^*-A^*\big)=0$.
Let $\big(Q(z)B^*-A^*\big)u=0$, $u\in \GG$. As
$\ker \big(\gr Q(z)-\Lambda^{A,B}\big)=0$, the corresponding
first component in $\big(B^*u,(Q(z)B^*-A^*\big)u\big)$
must vanish, i.e. $B^*u=0$, and then
$A^*u=0$. But $\ker A^*\cap\ker B^*=0$ due to~\eqref{eq-bg3},
and $u=0$. Therefore, the operator $Q(z)B^*-A^*$ has a bounded inverse due to
the closed graph theorem. Hence (a) implies (c).

Now let (c) hold, then we can rewrite~\eqref{eq-qab}
\begin{equation}
            \label{eq-loc3}
\gr
Q(z)-\Lambda^{A,B}=\Big\{\big(B^*(Q(z)B^*-A^*)^{-1}x,x\big),\,\,x\in\GG\Big\},
\end{equation}
and we get immediately $\ran \big(Q(z)-\Lambda^{A,B}\big)=\GG$
and $\ker \big(Q(z)-\Lambda^{A,B}\big)=0$, which exactly (a).
Therefore, we have shown that (a) is equivalent to (c).

Note that the condition (a) and, therefore, also (c), is invariant
under the change $z\leftrightarrow\bar z$, because they define the resolvent
set of the self-adjoint operator $H^{A,B}$, and the resolvent set
is symmetric under the complex conjugation.
The equivalence of (b) and (c) follows from the fact that
$Q(\bar z)B^*-A^*$ has $0$ in the resolvent set
if and only if its adjoint $BQ(z)-A$ has the same property
(here one can use the equality $Q(\overline z)=Q^*(z)$ following from~\eqref{q-fun}).
 
We have already proved the first equality in~\eqref{eq-qab}, see \eqref{eq-loc3}.
Let us show the second one. Let us the the notation of proposition~\ref{prop-dm}.
Replacing $z$ in \eqref{krein} by $\bar z$ and taking the adjoint on the both sides
one sees immediately that $C^*_\Lambda(\bar z)=C_\Lambda(z)$ for any $z\in\res H^0$.
On the other hand, we have shown already that for $\Lambda=\Lambda^{A,B}$
one has $C_\Lambda(z)=B^*\big(Q(z)B^*-A^*)^{-1}$, therefore,
$\big(BQ(z)-A\big)^{-1}B=C^*_\Lambda(\bar z)=C_\Lambda(z)=B^*\big(Q(z)B^*-A^*)^{-1}$.
\end{proof}

Now we are in position to reformulate proposition~\ref{prop-dm}
completely in the operator language, withour using linear relations.
\begin{thm}[Resolvent formula for normalized parameters]
       \label{prop-main}
Let $H^{A,B}$ be the self-adjoint extension of $S$ corresponding
to the boundary conditions~\eqref{A-B} with normalized $A$ and $B$.
A number $z\in\res H^0$ lies in $\spec H^{A,B}$
iff $0\in\spec \big(BQ(z)-A)$ or, equivalently, $0\in\spec \big(Q(z)B^*-A^*)$.
For any $z\in \res H^0\cap \res H^{A,B}$ there holds
\begin{align}
     \label{eq-krein1}
(H^{A,B}-z)^{-1}&=(H^0-z)^{-1}-\gamma(z)
B^*\big(Q(z)B^*-A^*\big)^{-1}\gamma^{\,*}(\bar{z}),\\
     \label{eq-krein2}
(H^{A,B}-z)^{-1}&=(H^0-z)^{-1}-\gamma(z)
\big(BQ(z)-A\big)^{-1}B\,\gamma^{\,*}(\bar{z}).
\end{align}
\end{thm}

\begin{proof}
Follows from lemma~\ref{lem-main}.
\end{proof}

Note that by setting $A=i(1+U)$, $B=1-U$ one obtains a global
expression for the resolvents which covers the whole family
of self-adjoint extensions.
A finite-dimensional case of this resolvent formula was obtained in~\cite{ap}
in the context of singular quantum-mechanical interactions.

\bigskip

Up to now, we used only normalized parameterizations for self-adjoint extensions.
As we see below, on example of elliptic boundary value problems,
normalized parameterization can be difficult to find. Let us formulate
the analogue of theorem~\ref{prop-main} for the case of non-normalized parameters.

\begin{thm}[Resolvent formula for non-normalized parameters]\label{prop-main2}
Let $H^{A,B}$ be the self-adjoint extension of $S$ corresponding
to some self-adjoint boundary conditions~\eqref{A-B}.
Then for any
$z\in\res H^{A,B}\cap \res H^0$ the operator $BQ(z)-A$ is
invertible on its domain and the equality~\eqref{eq-krein2} holds.
\end{thm}

\begin{proof}
Let $A_0$, $B_0$ be a strong parameterization of $\Lambda^{A,B}$. Due to proposition~\ref{prop-lab}
there exists a bounded injective operator $L$ such that
$A=L A_0$, $B=L B_0$. The resolvent formula~\eqref{eq-krein2} takes the form
\begin{align*}
R^{A,B}(z)\equiv R^{A_0,B_0}(z)&=R^0(z)-\gamma(z)\big(B_0 Q(z)-A_0\big)^{-1}B_0\gamma^*(\bar z)\\
&{}=R^0(z)-\gamma(z)\Big[L^{-1}\big(B Q(z)-A\big)\Big]^{-1}L^{-1} B\gamma^*(\bar z)\\
&{}=R^0(z)-\gamma(z)\big(BQ(z)-A\big)^{-1} B\gamma^*(\bar z).\qedhere
\end{align*}
\end{proof}

Although the resolvent formulas in theorems \ref{prop-main} and \ref{prop-main2}
have the same form, they are different from the point of view of the spectral analysis.
Namely, in the former case, the spectrum of $H^{A,B}$
in the gaps of $H^0$ is completely described of terms of the spectrum of $BQ(z)-A$.
In the latter case, the correspondence is more complicated. For example,
the operator $B Q(z)-A$ may have no bounded inverse for all $z$. Nevertheless,
one can describe at least the eigenvalues of $H^{A,B}$, cf. \cite[Theorem 1]{gm}
and \cite[Theorem 3.4]{pos2}.
\begin{thm}[Eigenvalues of self-adjoint extensions]
Let $H^{A,B}$ be the self-adjoint extension of $S$ corresponding to the boundary conditions \eqref{A-B}
with $A$ and~$B$ satisfying \eqref{eq-bg1} and~\eqref{eq-bg2}. The value $z\in\res H^0$
is an eigenvalue of $H^{A,B}$ iff $\ker \big(BQ(z)-A\big)\ne0$, and in this case one has
$\ker (H^{A,B}-z)=\gamma(z)\ker\big(BQ(z)-A\big)$.
\end{thm}
\begin{proof}
Let us show first that $\gamma(z)\ker\big(BQ(z)-A\big)\subset\ker (H^{A,B}-z)$.
For any $\xi\in\ker\big(BQ(z)-A\big)$ the element
$f(z,\xi):=\gamma(z)\,\xi$ is an eigenfunction of $S^*$ with the eigenvalue $z$,
because $\gamma(z)$ is an isomorphism between $\ker (S^*-z)$ and $\GG$.
Moreover, one has $\Gamma_1 f(z,\xi)=\Gamma_1 \Gamma_1^{-1}\xi=\xi$
and $\Gamma_2 f(z,\xi)=\Gamma_2\gamma(z)\xi=Q(z)\xi$, therefore,
$A\Gamma_1 f(z,\xi)-B\Gamma_2 f(z,\xi)=-\big(BQ(z)-A\big)\xi=0$, which means
that $f(z,\xi)$ is in the domain of $H^{A,B}$. Therefore, $f(z,\xi)$
is an eigenvector of $H^{A,B}$ with eigenvalue $z$.

Now let $z\in\res H^0$ be an eigenvalue of $H^{A,B}$ and
$f$ be a non-zero element of the corresponding subspace. Then $f$
is also an eigenvector of $S^*$. As $\gamma(z)$ is an ismorphism between $\ker (S^*-z)$
and $\GG$, there exists $\xi\in\GG\setminus \{0\}$ such that
$f=\gamma(z)\xi$. As previously, there holds $\Gamma_1 f=\xi$,
$\Gamma_2 f=Q(z)\xi$, and the condition $f\in \dom H^{A,B}$
takes the form $A\Gamma_1 f-B\Gamma_2f\equiv -\big(BQ(z)-A\big)\xi=0$.
\end{proof}

\section{Applications to quantum-mechanical solvable models}\label{sec4}

There are numerous works dedicated to the construction of boundary triples
in various situations, cf.~\cite{dm,dm2,gorb,pos1,pos2},
so we restrict ourselves by considering some situations where mixed
boundary conditions arise and the above results are useful.

\subsection{Generalized point interactions}
Krein's formula is used traditionally in the theory of zero-range potentials~\cite{alb}.
Let us illustrate the above construction by the generalized zero-range potentials.
Denote by $H^0$ the operator $-d^2/dx^2$ in $\HH=L^2(\RR)$.
Let $L$ be a uniformly discrete subset of $\RR$, i.e. there exists
$d>0$ such that $|a-b|\ge d$ for any $a,b\in L$ if $a\ne b$.
Denote by $S$ the restriciton of $H^0$ to the domain
$\dom S=\{f\in \dom H^0:\,f(a)=f'(a)=0\text{ for all }a\in L\}$.
Clearly, $S$ is a symmetric operator in $\HH$
and its deficiency indices are $(2|L|,2|L|)$ if $L$ is finite
and $(\infty,\infty)$ if $L$ is infinite, cf. \cite{Ku}.
Let us enumerate the points of $L$ in the increasing order,
$L=:\bigcup_{k\in K}\{a_k\}$, where $K$ is a subset of $\ZZ$
and $a_j< a_k$  if  $j<k$.

The adjoint operator $S^*$ is $-d^2/dx^2$ with the domain
$\dom S^*=W^{2,2}(\RR\setminus L)$.
The usual integration by parts and the Sobolev embedding theorem
show that the triple $(\GG,\Gamma_1,\Gamma_2)$
\begin{gather*}
\GG=l^2(K)\otimes\CC^2,\quad\Gamma_1f=\Big(f'(a_k-)-f'(a_k+),
f(a_k+)-f(a_k-)\Big)_{k\in K},\\
\Gamma_2f=\Big(\frac{f(a_k+)+f(a_k-)}{2},\frac{f'(a_k-)+f'(a_k+)}{2}\Big)_{k\in K},
\end{gather*}
is a possible choice of a boundary triple for $S$.
Note that the ``distinguished'' extension given by $\Gamma_1f=0$
is exactly $H^0$.

The general connecting self-adjoint boundary conditions at the points of $L$ have the form~\cite{ADK,PS}
\begin{gather*}
\begin{pmatrix}
f'(a_k+)\\f(a_k+)
\end{pmatrix}=
e^{i\theta_k}
\begin{pmatrix}
\alpha_k & \beta_k\\
\gamma_k & \delta_k
\end{pmatrix}
\begin{pmatrix}
f'(a_k-)\\f(a_k-)
\end{pmatrix},\\
\theta_k\in[0,\pi],
\quad \alpha_k,\beta_k,\gamma_k,\delta_k\in\RR,
\quad \alpha_k\delta_k-\beta_k\gamma_k=1,\quad k\in K,
\end{gather*}
and can be rewritten as $A\Gamma_1 f=B\Gamma_2 f$, where $A$ and $B$ are block-diagonal matrices, $A=\diag(A_k)$, $B=\diag (B_k)$,
\[
A_k=\begin{pmatrix}
1+\alpha e^{i\theta_k} & \beta_k e^{i\theta_k}\\
\gamma e^{i\theta_k} & \delta_k e^{i\theta_k}-1
    \end{pmatrix},\quad
B_k=\begin{pmatrix}
\beta_k e^{i\theta_k} & 1-\alpha_k e^{i\theta_k}\\
1+\delta_k e^{i\theta_k} & -\gamma_k e^{i\theta_k}
    \end{pmatrix}.
\]
Denote by $H$ the resulting self-adjoint operator.

Let $h=(h',h'')$, $h',h''\in l^2(K)$.
By direct calculations one shows that the solution of $(S^*-z)g=0$ with $\Gamma_1 g=h=(h',h'')$
and $z\notin\spec H^0\equiv [0,+\infty)$ has the form
\[
g(x)=\gamma(z)h(x)=\frac{1}{2}\,\sum_{k\in K} \Big(
\dfrac{h'_k}{\sqrt{-z}}+h''_k \sgn (x-a_k)\Big)\,e^{-\sqrt{-z}|x-a_k|}
\]
For the corresponding Krein matrix one has $Q(z)=\big(Q^{jk}(z)\big)_{j,k\in K}$,
where the $2\times 2$ blocks $Q^{jk}(z)$ are:
$Q^{jk}(\zeta)$,
\[
Q^{jk}(z)=\dfrac{e^{-\sqrt{-z}|a_j-a_k|}}{2}\,
\begin{pmatrix}
\dfrac{1}{\sqrt{-z}} & \sgn (a_j-a_k)\\[\bigskipamount]
-\sgn(a_j-a_k) & -\sqrt{-z}\,
\end{pmatrix}. 
\]
The resolvent formula presented above can handle arbitrary boundary conditions
at the points of $L$. In particular, if $L$ is finite, the equation
$\det\big[BQ(z)-A\big]=0$ provides a condition for $z$ to be an eigenvalue of $H$.
A part of this construction can be transferred to matrix-valued problems
with point interactions~\cite{KP}.

\subsection{Hybrid manifolds}

Another object on which mixed boundary conditions naturally arise
is delivered by Schr\"odinger operators acting on a space consisting
of pieces with different dimensions. We consider an example of a space consisting
just of two pieces of different dimensions using the construction of boundary triple
from \cite{bg,E05} using some asymptotic properties of the Green function from \cite{bgp}.
Let $M$ and $N$ be two- and three-dimensional
complete Riemannian manifolds of bounded geometry, respectively. (The class of manifolds of bounded geometry
includes compact manifolds, homogeneous spaces with invariant metrics, etc., see \cite{bgp}).
Let $V\in L^2_\text{loc}(M)$, $W\in L^2_\text{loc}(N)$
be semibounded below. The expressions $H_M=-\Delta_M+V$ and $H_N=-\Delta_N+W$,
where $\Delta_{M}$  and $\Delta_{N}$ are the Laplace-Beltrami operator on $X$ and $Y$, respectively,
define Schr\"odinger operators acting on $L^2(N)$ and $L^2(N)$.
By identifying two marked
points $m\in M$ and $n\in N$ one arrives at a new topological space.
Our aim to define Schr\"odinger operator on this new space using the partial operators $H_{M}$, $H_N$.
The new phase space is $\HH=L^2(M)\oplus L^2(N)$. Denote by $S_{M/N}$ the restrictions
of $H_{M/N}$ to functions vanishing at $m/n$. These operators are symmetric and have deficiency indices
$(1,1)$. By definition, Schr\"odinger operators on $\HH$ are self-adjoint extension
of $S=S_N\oplus S_M$. The corresponding boundary triples can be constructed as follows.

Let $G_{M/N}(x,y;z)$ be the Green function of $H_{M/N}$, i.e. the integral kernel
of the corresponding resolvents $(H_{M/N})^{-1}$ for $z\notin\spec H_{M/N}$. These kernels
are continuous for $x\ne y$, and there exist functions $F_{M/N}(x,y)$
such that $G_{M/N}(x,y;z)=F_{M/N}(x,y)+G^\text{ren}_{M/N}(x,y;z)$, where the renormalized
Green function $G^\text{ren}_{M/N}(x,y;z)$ is continuous in the whole space $M\times M$
or $N\times N$, see \cite{bgp}. Moreover, one can always put $F_M(x,y)=-\dfrac{1}{2\pi}\log d(x,y)$
and, if $W\in L^p_\text{loc}(N)$ with $p>3$, $F_{N}(x,y)=\dfrac{1}{4\pi d(x,y)}$.
The domain of the adjoint  operators are $\dom S^*_{M}=\dom H_M+\CC G_M(\cdot,m;z)$
and $\dom S^*_{N}=\dom H_N+\CC G_N(\cdot,n;z)$ with any $z\notin\spec H_{M/N}$.
For any $f\in \dom S^*_{M/N}$ one has the following asymptotic
expansion near $m/n$: $f(x)=a_{M/N}(f)F_{M/N}(x,m/n)+b_{M/N}(f)+o(1)$, $a_{M/N}$, $b_{M/N}$
are linear functionals on $\dom S^*_{M/N}$. Using the partial integration
one can see that $(\CC,\Gamma^{M/N}_1,\Gamma^{M/N}_2)$ with $\Gamma_1^{M/N}f=a_{M/N}(f)$,
$\Gamma_1^{M/N}f=b_{M/N}(f)$ are boundary triples for $S_{M/N}$, cf. lemma~5 in \cite{bg}.
The corresponding $\Gamma$-fields $\gamma_{M/N}$ and $\mathcal Q$-functions $Q^{M/N}$ are
$\gamma_{M/N}(z):\CC\ni\xi\mapsto G_{M/N}(\cdot,m/n;z)$ and $Q^{M/N}(z)=G^\text{ren}(m/n,m/n;z)$. 
Clearly, $(\CC^2,\Gamma_1,\Gamma_2)$ with $\Gamma_j=\Gamma_j^M\oplus \Gamma_j^N$, $j=1,2$,
is a boundary triple for $S$, and the corresponding $\Gamma$-field and the $\mathcal Q$-function
are the direct sums of those for $S_{M/N}$. The most general boundary conditions
are $A\Gamma_1 f=B\Gamma_2$ with $2\times 2$ matrices $A$ and $B$ with $AB^*=BA^*$
and $\det (AA^*+BB^*)\ne 0$. Using the resolvent formulas~\ref{eq-krein1}
or \eqref{eq-krein2} one can perform a complete spectral analysis of the operators obtained
by repeating the arguments of \cite{bg} where special boundary conditions were studied.
The construction presented admits a natural generalization
to infinite direct sums and periodic hybrid structures, cf. \cite{bgl}.

\subsection{Laplacian on the half-space}

One of related examples comes from the Robin-type elliptic boundary value problems,
i.e. the Laplacians with boundary conditions of the form 
$au|_\Gamma+b\dfrac{\partial u}{\partial n}\big|_\Gamma=0,$
where $\Gamma$ is a boundary of a certain
domain $\Omega$ and $a$, $b$ are real-valued coefficients.
If the coefficient $b$ vanishes on some set of non-zero measure,
one has the so-called mixed boundary problem; such boundary conditions
arise in various areas of geometric analysis
and mathematical physics~\cite{AIM,BE,DK,dow}. We discuss here an explicitly solvable case,
when $\Omega$ is a half-space. Similar constructions can be done
also done in tube-like domains, see section~9 in~\cite{dm}
and section~\S4.7 in~\cite{gorb} for discussion. 

Denote $\RR^n_+:=\{x=(x_1,\dots,x_n)\in\RR^n:\, x_n\ge 0\}$, $n\ge 2$.
By $S$ denote the Laplace operator
with the domain $\dom S=\{f\in H^2(\RR^n_+):\, f|_{x_n=0}=\partial f /\partial x_n\big|_{x_n=0}=0\}$.
Clearly, $S$ is a symmetric operator in $\HH:=L^2(\RR^n_+)$ with infinite deficiency indices.

Let us consider the self-adjoint extensions of $S$ corresponding to the Robin-type boundary conditions,
\begin{equation}
      \label{eq-rob2}
af|_{x_n=0}+b\,\frac{\partial f}{\partial x_n}\big|_{x_n=0}=0,\quad a.e.
\end{equation}
We assume that $a\in C^1(\RR^{n-1})$, $b\in C^2(\RR^{n-1})$ are real-valued function,
and that $a$, $a'$, $b$, $b'$, $b''$ are bounded. Both $a$ and $b$ can vanish
on sets of non-zero measure, but the sum $|a|+|b|$ must be nonzero almost everywhere.

Denote by $H^{a,b}$ the self-adjoint operator introduced above and
by $H^0$ the Dirichlet Laplacian. Our aim is to express the resolvent
of $H^{a,b}$ in terms of $H^0$, $a$ and $b$.
We use the construction of section~9 in~\cite{dm} for the boundary
triple and the Krein functions.

For any $z<0$ denote $\Lambda(z)=(-\Delta_{n-1}-z)^{1/2}$, where $\Delta_{n-1}$
in the Laplacian in $\RR^{n-1}$. Clearly,
$\Lambda(z)$ is an injective non-negative self-adjoint operator in $L^2(\RR^{n-1})$.
Fix any $\lambda<0$, then as a boundary triple $(\GG,\Gamma_1,\Gamma_2)$ for $S$ one can take
\begin{multline}
             \label{eq-hpgg}
\GG=L^2(\RR),\quad
(\Gamma_1 f)(x)= \big(\Lambda^{-1/2}(\lambda)\big)_x f|_{x_n=0},\\
(\Gamma_2 f)(x)=\big(\Lambda^{1/2}(\lambda)\big)_x \Big(
\frac{\partial f}{\partial x_n}\Big|_{x_n=0}+\Lambda(\lambda) f|_{x_n=0}\Big).
\end{multline}
The corresponding map $\gamma(z)$ represent the well-known formula
for the solution of the Dirichlet problem,
\[
\gamma(z) u=-\int_{\RR^{n-1}}\frac{\partial G (\cdot,y;z)}{\partial y_n}\big|_{y_n=0}
\Lambda^{1/2}(\lambda)u(y)\,dy,\quad u\in L^2(\RR^{n-1}),
\]
where $G(x,y;z)$ is the Green function for $H^0$.
The corresponding $Q$-function takes the form
$Q(z)=\big(\Lambda(\lambda)-\Lambda(z)\big)\Lambda(\lambda)$.
Substituting \eqref{eq-hpgg} into \eqref{eq-rob2}
applying on the both sides the operator $\Lambda^{-3/2}(\lambda)$
we obtain the boundary conditions $A\Gamma_1 f=B\Gamma_2 f$ with
\begin{equation}
              \label{eq-hpAB}
A=\Lambda^{-3/2}(\lambda)\big(b\Lambda(\lambda)-a\big)\Lambda^{1/2}(\lambda),
\quad B=\Lambda^{-3/2}(\lambda)b\Lambda^{-1/2}(\lambda).
\end{equation}
As the multiplications by $a$ and $b$ are bounded operators in $H^{-1}(\RR)$
and $H^{-2}(\RR)$, respectively, the coefficients $A$ and $B$ in \eqref{eq-hpAB}
are bounded operators in $L^2(\RR)$.
In general, $A$ and $B$ are not normalized
and it is very difficult to normalize them in a closed form.
Therefore, we can obtain only the weak form of the resolvent formula~\eqref{eq-krein2},
see theorem~\ref{prop-main2}.
The expression $B Q(z)-A$ entering~\eqref{eq-krein2} takes the form
$B Q(z)-A=\Lambda^{-3/2}(\lambda)\big(a-b\Lambda(z)\big)\Lambda^{1/2}(\lambda)$
In particular, a real number $z<0$ is an eigenvalue
of $H^{a,b}$ iff there exists $u\in H^{-1/2}(\RR)$ with
$\big(a-b\Lambda(z)\big) u=0$.

\section*{Acknowledgments} The work was supported by the Deutsche Forschungsgemeinschaft,
the Sonderforschungsbereich ``Raum, Zeit, Materie'' (SFB 647, Berlin),
and the International Bureau of BMBF at the German Aerospace Center
(IB DLR, cooperation Germany -- New Zealand NZL 05/001).

\end{document}